\documentclass[12pt,a4paper,reqno]{amsart}
\usepackage{amsmath}
\usepackage{amsfonts}
\usepackage{amssymb}
\usepackage{bm} 

\newcommand\eps{\varepsilon}
\newcommand\R{{\mathbf{R}}}
\newcommand\C{{\mathbf{C}}}
\newcommand\Z{{\mathbf{Z}}}

\newcommand\N{{\mathcal{N}}}
\renewcommand\S{{\mathcal{S}}}
\newcommand\sgn{{\operatorname{sgn}}}

\theoremstyle{plain}
  \newtheorem{theorem}[subsection]{Theorem}
  
  \newtheorem{proposition}[subsection]{Proposition}
  \newtheorem{lemma}[subsection]{Lemma}
  \newtheorem{corollary}[subsection]{Corollary}

\theoremstyle{remark}
  \newtheorem{remark}[subsection]{Remark}

\theoremstyle{definition}
  \newtheorem{definition}[subsection]{Definition}

\include{psfig}

\begin{document}

\title[Scattering of quartic gKdV equation]{Scattering for the quartic generalised Korteweg-de Vries equation}
\author{Terence Tao}
\address{Department of Mathematics, UCLA, Los Angeles CA 90095-1555}
\email{tao@math.ucla.edu}

\vspace{-0.3in}
\begin{abstract}
We show that the quartic generalised KdV equation
$$ u_t + u_{xxx} + ( u^4 )_x = 0$$
is globally wellposed for data in the critical (scale-invariant) space $\dot H^{-1/6}_x(\R)$ with small norm (and locally wellposed for large norm), improving a result of Gruenrock \cite{gruenrock}.  As an application we obtain scattering results in $H^1_x(\R) \cap \dot H^{-1/6}_x(\R)$ for the radiation component of a perturbed soliton for this equation, improving the asymptotic stability results of Martel and Merle \cite{martel}.
\end{abstract}

\maketitle

\section{Introduction}

In this paper we study the global behaviour for the quartic generalised KdV (gKdV) Cauchy problem
\begin{equation}\label{eq:gkdv}
 u_t + u_{xxx} + ( u^4 )_x = 0; \quad u(0,x) = u_0(x) \in H^s_x(\R \to \C),
\end{equation}
where $u: \R \times \R \to \C$ is a complex-valued function.  Note that as the power of the nonlinearity is even, this equation should not be considered either focusing or defocusing, although it does admit soliton solutions; in particular, the sign of the nonlinearity $(u^4)_x$ is not relevant as can be seen by the transformation $u \mapsto -u$.

The Cauchy problem \eqref{eq:gkdv} was shown to be locally wellposed
in $H^s_x(\R \to \C)$ for\footnote{At first glance it is not immediately obvious that the equation should even make sense distributionally for negative regularities.  However, the smoothing effect inherent in the Strichartz estimates will allow us to (for instance) place $u$ in $L^6_{t,x}$ whenever $u_0 \in \dot H^{-1/6}_x$.} $s > -1/6$ and globally wellposed in $H^s(\R \to \R)$ for $s \geq 0$ by Gruenrock \cite{gruenrock} (building on an earlier local and global results for $s \geq \frac{1}{12}$ and $s \geq 1$ respectively by Kenig, Ponce, and Vega \cite{kpv}).  
The exponent $s=-1/6$ is \emph{critical} for this equation, because the homogeneous Sobolev space $\dot H^{-1/6}_x(\R \to \C)$
is invariant under the scaling symmetry
\begin{equation}\label{scaling}
u(t,x) \mapsto \frac{1}{\lambda^{2/3}} u(\frac{t}{\lambda^3}, \frac{x}{\lambda}); \quad u_0(x) \mapsto \frac{1}{\lambda^{2/3}} u_0(\frac{x}{\lambda})
\end{equation}
of the equation \eqref{eq:gkdv}.  

The above local wellposedness results were obtained by the contraction mapping method.  By a refinement of the spaces and estimates used in \cite{gruenrock}, we are able to also establish a well-posedness result at the critical regularity:

\begin{theorem}[Critical wellposedness]\label{main} The equation \eqref{eq:gkdv} is locally wellposed in $\dot H^{-1/6}_x(\R \to \C)$, with the time of existence being infinite when the $\dot H^{-1/6}_x(\R \to \C)$ norm is sufficiently small.  
\end{theorem}

We shall define precisely what we mean by ``wellposedness'' in Theorem \ref{cle} below.  As one might expect from a critical regularity iteration argument, our methods also yield persistence of regularity, as well as and scattering when the critical norm $\|u_0\|_{\dot H^{-1/6}_x(\R \to \C)}$ is small; on the other hand, for large data the time of existence depends on the frequency profile of the data $u_0$, and not just on the norm (again, see Theorem \ref{cle} for a precise statement).  As our methods rely on a pure iteration argument, the solution will depend analytically on the initial data, and one can extend the results without difficulty to more general systems with $(u^4)_x$ type nonlinearity where $u$ is vector-valued.

Our iteration methods rely on a combination of linear, bilinear, and quartilinear Strichartz estimates, together with some
basic critical $X^{s,b}$ (or more precisely, $\dot X^{s,b,q}$) theory.  Interestingly, for this wellposedness result we do not exploit local smoothing or maximal function estimates, although the Strichartz estimates that we use do capture some smoothing effects of the equation, and when we incorporate the influence of a soliton (see below) then the Kato local smoothing effect will become crucial.

From a numerology viewpoint, this appears to be the first critical \emph{negative-regularity} wellposedness result for a dispersive equation (excluding artificial examples, such as taking an $L^2_x$-critical equation and conjugating it by an inverse derivative to make it $\dot H^{-1}_x$-critical).  The point seems to be that the usual obstructions to reaching a negative critical regularity, such as Galilean invariance or (more generally) the self-interaction of high frequencies, do not seem to be significant for this equation, in large part because the power of the quartic nonlinearity $(u^4)_x$ is \emph{even}\footnote{The original Korteweg de Vries (KdV) equation, with a nonlinearity of $(u^2)_x$ type, has similarly weak self-interactions of high frequency when expanded to the first nonlinear iterate, but when one then computes the second nonlinear iterate one sees significant self-interaction again, which ultimately blocks iteration methods at $s=-3/4$ rather than the critical $s=-3/2$ (see \cite{cct}).  The situation here can be clarified by inverting the Miura transform to convert the KdV equation to the modified KdV (mKdV) equation, with nonlinearity of $(u^3)_x$ type.  In the quartic case, the higher critical regularity $s = -1/6$ and the higher power of the nonlinearity seems to make this second iteration self-interaction effect insignificant.}.

Now we combine this perturbative analysis with the more global analysis of Martel and Merle (based largely on conservation laws, virial identities, and linearised stability analysis, rather than on iteration in carefully chosen function spaces), specialising now to the finite energy real-valued case.  More precisely, we consider small $H^1_x(\R \to \R)$ perturbations of the soliton solution
$$ u(t,x) := Q(x - t)$$
where 
$$ Q(x) := \left( \frac{4}{2 \cosh^2(\frac{3}{2}x)} \right)^{1/3}$$
is the unique positive even Schwartz solution to the ODE
\begin{equation}\label{soliton-eq}
Q'' + Q^4 = Q.
\end{equation}

The following result was established by Martel and Merle \cite{martel}, \cite{martel2}:

\begin{theorem}[Asymptotic stability of solitons]\label{asymptotic} Suppose that $u_0 \in H^1_x(\R \to \R)$ is such that $\| u_0 - Q \|_{H^1_x(\R)} \leq \eps$ for some sufficiently small $\eps > 0$.  Then the global solution $u$ to \eqref{eq:gkdv} admits a decomposition
$$ u(t,x) = \frac{1}{\lambda^{2/3}(t)} Q\left( \frac{x-x(t)}{\lambda(t)} \right) + w(t,x)$$
where $\lambda: \R \to \R^+$, $x: \R \to \R^+$ are $C^1$ functions such that
\begin{equation}\label{lambda}
 \lambda(t) = 1 + O(\eps); \quad x'(t) = 1 + O(\eps) \hbox{ for all } t
 \end{equation}
and $w$ obeys the estimates
\begin{equation}\label{lh1}
\| w \|_{L^\infty_t H^1_x(\R)} \lesssim \eps
\end{equation}
and
\begin{equation}\label{wox}
 \int_\R \int_\R (|w(t,x)|^2 + |w_x(t,x)|^2) e^{-\sigma |x-x(t)|}\ dx dt \lesssim_\sigma \eps^2
\end{equation}
for any $\sigma > 0$.  Furthermore, we have the estimate
\begin{equation}\label{lambda-rate}
|\lambda'(t)| + |x'(t) - \lambda(t)^{-2}|^2 \lesssim
\int_\R |w(t,x)|^2 e^{-\frac{1}{2} |x-x(t)|}\ dx 
\end{equation}
\end{theorem}

\begin{proof} See \cite[Lemma 1]{martel2} and \cite[(3.2)]{martel2}. 
\end{proof}

\begin{remark}\label{mm} In \cite{martel}, \cite{martel2} some additional information is deduced on $\lambda(t), x(t), w(t)$; for instance, the function $\lambda(t)$ (which measures the spatial width of the soliton) converges to a limit $\lambda(\pm \infty)$ as $t \to +\infty$ or as $t \to -\infty$.  Similarly, the quantity $x'(t)$ (which measures the velocity of the soliton) also converges to $\lambda(\pm \infty)^{-2}$ as $t \to \pm \infty$.  On the other hand, one does not expect $x(t) - \lambda(\pm \infty)^{-2} t$ to stay bounded for this type of (slowly decaying) perturbation, see \cite{martel2}.  Also, as $t \to +\infty$, $w(t)$ converges weakly in $H^1_x(\R \to \R)$ to zero and strongly in
$H^1_x([\beta t, +\infty) \to \R)$ for some $\beta = \beta(\eps) > 0$ which goes to zero as $\eps \to 0$; similarly as $t \to -\infty$.  The results are not specific to the quartic gKdV and also hold for the other $L^2$-subcritical equations, namely the quadratic (KdV) and cubic (mKdV) equations.
\end{remark}

Using our estimates, we are able to obtain some further control on $w$, showing that it scatters to a free solution:

\begin{theorem}[Scattering]\label{main2} Let the notation and assumptions be as in Theorem \ref{asymptotic}.  Assume also that $\|u_0 - Q \|_{\dot H^{-1/6}_x(\R \to \R)} \lesssim \eps$.
Then there exists $w_+ \in H^1_x(\R \to \R) \cap \dot H^{-1/6}_x(\R \to \R)$ with 
$$\|w_+\|_{H^1_x(\R \to \R) \cap \dot H^{-1/6}_x(\R \to \R)} = O(\eps)$$ 
such that
$$ \lim_{t \to +\infty} \| w(t) - e^{-t\partial_{xxx}} w_+ \|_{H^1_x(\R \to \R) \cap \dot H^{-1/6}_x(\R \to \R)} = 0.$$
\end{theorem}

\begin{remark}
It will be clear from the proof that $w$ will also obey all the standard estimates that a linear $H^1 \cap \dot H^{-1/6}$ solution to the Airy equation $w_t + w_{xxx} = 0$ with norm $O(\eps)$ would obey, such as Strichartz estimates, local smoothing estimates, maximal function estimates, etc.  For instance, we will have
the global spacetime integrability estimate
$$ \| w \|_{L^6_{t,x}(\R \times \R)} \lesssim \eps.$$
The additional condition $\|u_0 - Q \|_{\dot H^{-1/6}_x} \lesssim \eps$ is a decay condition on the very low frequencies of $u_0$, and is necessary to obtain the type of global asymptotic results here because one must have $w$ small in at least one scale-invariant norm in order to have any hope of closing a global iteration argument.  Note that from Sobolev embedding it would be enough for $u_0$ to be close to $Q$ in $L^{3/2}_x$ (and in $H^1_x$, of course).
\end{remark}

\begin{remark} This result complements the results of C\^ote \cite{cote}, \cite{cote2} constructing large data wave operators for the equation \eqref{eq:gkdv}, both with and without the presence of one or more solitons.  In those works, the scattering state $w_+$ was assumed to have some additional decay at spatial infinity (e.g. $H^{1,1}_x$); it now seems likely that with the methods here, these decay conditions could be removed, although we do not pursue this issue here.  See also \cite{cote3}, which (similarly to this paper) also combines the stability analysis of solitons with a critical wellposedness theory.
\end{remark}

\begin{remark} As $t \to +\infty$, the analysis in \cite{martel}, \cite{martel2} gives satisfactory control of the radiation $w$ in the region $x > \beta t$ (which is where the soliton lies); Theorem \ref{main2} complements this control by showing in the complementary region $x \leq \beta t$ that the solution behaves like a linear solution to the Airy equation.  In particular it is now not hard to show that we can take $\beta = 0$ in Remark \ref{mm}.  
\end{remark}

\begin{remark} Our methods exploit the even power of the nonlinearity $(u^4)_x$, which ensures that the self-interaction of any given frequency mode is either highly non-resonant or low frequency, and thus negligible in both cases.  This lack of strong self-interaction allows us to exploit bilinear Strichartz estimates effectively.  Thus the methods here do not yield any immediate progress on the modified KdV (mKdV) equation (which is \eqref{eq:gkdv} but with nonlinearity $\pm (u^3)_x$), in which local wellposedness is only known in $H^s_x$ for $s \geq 1/4$ despite the presence of a subcritical $L^2_x$ conservation law.  Indeed it is known for that equation that a pure iteration method cannot work for 
any $s < 1/4$, precisely because of the self-interaction of frequency modes; see \cite{cct}.  There is however some hope that one could renormalise away this self-interaction to break the $s=1/4$ barrier; see \cite{hideo} for an execution of this idea in the periodic setting.
\end{remark}

\begin{remark} The conclusion of Theorem \ref{main2} can be written in the form
$$ u(t,x) = \frac{1}{\lambda^{2/3}(t)} Q\left( \frac{x-x(t)}{\lambda(t)} \right) + e^{-t\partial_{xxx}} w_+(x) 
+ o_{t \to +\infty}(1)_{H^1_x}$$
where we use $o_{t \to +\infty}(1)_{H^1_x}$ to denote a quantity which goes to zero in $H^1_x$ norm
as $t \to \infty$.  One can then use the conservation of mass
$$ \int_\R u(t,x)^2\ dx$$
and energy
$$ \int_\R \frac{1}{2} u_x(t,x)^2 - \frac{1}{5} u(t,x)^5\ dx$$
to obtain asymptotic decoupling of mass and energy:
\begin{align*}
\int_\R u_0(x)^2\ dx &= \lambda(+\infty)^{-1/3} \int_\R Q(x)^2\ dx + \int_\R w_+(x)^2\ dx \\
\int_\R \frac{1}{2} u'_0(x)^2 - \frac{1}{5} u_0(x)^5\ dx &= 
\lambda(+\infty)^{-7/3} \frac{1}{10} \int_\R Q(x)^2\ dx + \int_\R w'_+(x)^2\ dx;
\end{align*}
we leave the details to the reader.  
\end{remark}

We prove Theorem \ref{main2} in Section \ref{main2-sec}, assuming some technical estimates which we then prove in Section \ref{estimates-sec}.  The idea is to subtract off the soliton component $\frac{1}{\lambda^{2/3}(t)} Q( \frac{x-x(t)}{\lambda(t)} )$ and analyse the radiation component $w$ directly via the global perturbation methods used to
establish Theorem \ref{main}.  A difficulty arises because the interaction of the soliton and the radiation (and also the radiation caused by the soliton changing in size or having unexpected velocity) decays in time only in an $L^2_t$ sense rather
than an $L^1_t$ sense (thanks to \eqref{wox}); also, the regularity of this interaction forcing term is one derivative worse than expected ($L^2_x$ rather than $H^1_x$).  However, the crucial fact that the soliton moves to the \emph{right} (while the fundamental solution for the Airy equation essentially moves to the \emph{left}) will yield both a local smoothing effect and a certain almost orthogonality
which will allow us to recover the derivative and sum the $L^2_t$ forcing term\footnote{One way to view this is by observing that the interaction between the soliton and radiation, while not in a traditional energy space such as $L^1_t L^2_x$, is in a nontraditional energy space $L^1_h L^2_t$, where we use the curvilinear coordinate system $(h,t) := (x - x(t),t)$.  Such ``tilted energy spaces'' seem to appear frequently in critical wellposedness theory in which derivatives are present in the nonlinearity; see also \cite{tataru:wave2}, \cite{tao:wavemap2}, \cite{ik}, \cite{bej} for further examples. The situation here is slightly simpler than in those papers because we do not really encounter angular separation issues in one dimension, except via our use of the Riesz transforms $P_-$, $P_+$ to separate positive and negative frequencies.}.  This however forces a certain technical expansion of the
function spaces used to iterate in, making the argument more complicated than that in Theorem \ref{main}.

\section{Notation}\label{notation-sec}

We may take all functions to be smooth and rapidly decreasing in space to facilitate the rigorous justification of
various steps in the argument; once one establishes uniform estimates in this case, one can then pass to rough solutions
by the usual limiting argument.  As it is the estimates which are at the heart of the analysis, we shall thus gloss over
the justification of various steps (such as interchanging two integrals) by implicitly assuming as much regularity and
spatial decay as necessary\footnote{While we do not have infinite decay of solutions in time, in practice we can circumvent this by restricting the time axis $\R$ to some large compact interval $[-T,T]$ (possibly using smooth cutoffs if desired), establishing estimates which are asymptotically uniform in $T$, and then letting $T$ go to infinity.  As these technical steps are rather standard and uninteresting we will not pursue them explicitly here.}.  

We use $X \lesssim Y$ or $X = O(Y)$ to denote the estimate $X \leq C Y$ for some constant $C > 0$.  Occasionally our constants shall depend on an additional parameter such as $s$, in which case we subscript the $\lesssim$ or $O()$ notation accordingly, thus for instance $X \lesssim_s Y$ denotes the estimate $X \leq C(s) Y$ for some $C(s) > 0$ depending only on $s$.

We use subscripting by $t$ and $x$ to denote partial differentiation, and primes to denote ordinary differentiation.

We normalise our spatial and spacetime Fourier transforms
$$ \hat f(\xi) := \frac{1}{2\pi} \int_\R e^{-i x \xi} f(x)\ dx$$
and
$$ \tilde u(\tau,\xi) := \frac{1}{(2\pi)^2} \int_\R \int_\R e^{-i (x \xi + t \tau)} u(t,x)\ dx dt$$
in order to obtain the inversion formulae
$$ f(x) = \int_\R e^{ix \xi} \hat f(\xi)\ d\xi$$
and
$$ u(t,x) = \int_\R \int_\R e^{i(x \xi+t\tau)} \tilde u(\tau, \xi)\ d\xi.$$
We use the usual Lebesgue norms
$$ \| f \|_{L^p_x(\R \to \C)} := (\int_\R |f(x)|^p\ dx)^{1/p}$$
and
$$ \| u \|_{L^q_t L^r_x(\R \times \R \to \C)} := (\int_\R \|u(t)\|_{L^r_x(\R \to \C)}^q)^{1/q}$$
with the usual modifications when $p,q,r = \infty$.  We shall omit the domain and range of these spaces when they are clear from context.  We also abbreviate $L^p_t L^p_x$ as $L^p_{t,x}$.

We use the Japanese bracket $\langle x \rangle := (1 + |x|^2)^{1/2}$, and use this to define the homogeneous and inhomogeneous differential operators $|\nabla|^s$, $\langle \nabla \rangle^s$ via the Fourier transform as
$$ \widehat{|\nabla|^s f}(\xi) := |\xi|^s \hat f(\xi); \quad \widehat{\langle \nabla\rangle^s f}(\xi) := \langle \xi \rangle^s \hat f(\xi).$$
We then define the Sobolev spaces $\dot H^s_x = \dot H^s_x(\R \to \C)$, $H^s_x = H^s_x(\R \to \C)$ as
$$ \| f \|_{\dot H^s_x(\R \to \C)} := \| |\nabla|^s f \|_{L^2_x}; \quad
\| f \|_{H^s_x(\R \to \C)} := \| \langle \nabla \rangle^s f \|_{L^2_x}.$$
We also need the Riesz transforms $P_-, P_+$ defined as
$$ \widehat{P_- f}(\xi) = 1_{\xi < 0} \hat f(\xi); \quad \widehat{P_+ f}(\xi) = 1_{\xi \geq 0} \hat f(\xi)$$
and the propagators
$$ \widehat{e^{-t\partial_{xxx}} u_0}(\xi) := e^{it\xi^3} \hat u_0(\xi).$$
We observe the \emph{Duhamel formula}
$$ u(t) = e^{-t\partial_{xxx}} u(0) + \int_0^t e^{-(t-t')\partial_{xxx}} (\partial_t + \partial_{xxx}) u(t')\ dt'$$
where we adopt the convention that $\int_0^t = -\int_t^0$ when $t < 0$.

For minor technical reasons we shall use a slightly unusual Littlewood-Paley decomposition, using powers of $1.001$ instead of $2$.  Let $\varphi: \R \to \R$ be a smooth even function with $\varphi(\xi) = 1$ for $\xi \leq 1$ and $\varphi(\xi) = 0$ for $|\xi| > 1.001$.  Whenever $N$ is a power of $1.001$, we define the operator
$$ \widehat{P_{\leq N} f}(\xi) := \varphi(\frac{\xi}{N}) \hat f(\xi)$$
and then
$$ P_N f:= P_{\leq N}f - P_{\leq N/1.001} f; \quad P_{>N} f := f - P_{\leq N} f.$$
We remark that these operators are convolutions with absolutely integrable kernels and are thus bounded on every
translation-invariant Banach space, uniformly in $N$.  Similarly $\partial_x P_N$ is bounded on every such space with an 
operator norm of $O(N)$.  We also observe the Bernstein inequality
$$ \| P_N f\|_{L^q_x} \lesssim_{p,q} N^{1/p - 1/q} \|P_N f\|_{L^p_x}$$
whenever $1 \leq p \leq q \leq \infty$.

Henceforth $N$ is always understood to be a power of $1.001$.

For $b \in \R$ and $1 \leq q \leq \infty$, we define the $\dot X^{0,b,q} = \dot X^{0,b,q}(\R \times \R \to \C)$ norm of a function $u$ via its spacetime Fourier transform as
$$ \| u \|_{\dot X^{0,b,q}} := (\sum_{k=-\infty}^{+\infty} 
(2^{bk} \| \tilde u \|_{L^2_{\tau,\xi}(A_k)})^q)^{1/q}$$
where $A_k$ is the region 
$$ A_k := \{ (\tau,\xi):2^k \leq |\tau-\xi^3| < 2^{k+1} \}.$$
We define $\dot X^{0,b,q}$ to be the weak closure of the test functions that are uniformly bounded in the above norm.  In particular, when $b = 1/2$ and $q = +\infty$, we see that $\dot X^{0,1/2,\infty}$ contains the free solutions $e^{-t\partial_{xxx}} u_0$, where $u_0$

If $X$ and $Y$ are two function spaces on the same domain, we use $X \cap Y$ to denote the function space of functions in both $X$ and $Y$ with norm $\| f \|_{X \cap Y} := \|f\|_X + \|f\|_Y$, and $X+Y$ to denote the function space of sums of functions in $X$ and functions in $Y$ with norm
$$ \|f\|_{X+Y} := \inf \{ \|f_1\|_X + \|f_2\|_Y: f = f_1 + f_2 \}.$$
Also, if $\Omega$ is a subdomain, we use $X(\Omega)$ to denote the restriction of the functions in $X$ to $\Omega$, with norm
$$ \|f\|_{X(\Omega)} := \inf \{ \|g\|_X: f = g|_\Omega \}.$$

\section{Free evolution estimates}

In this section (and the next two) all functions are allowed to be complex-valued.

We now record some standard (and less standard) estimates for the free propagator $e^{-t\partial_{xxx}}$ for the Airy equation.
For technical reasons (relating to the Fourier transform) we shall need to consider complex-valued solutions, even though all
our applications are for real-valued functions.  We make the trivial remark that if $u(t) = e^{-t\partial_{xxx}} u_0$, then
$\overline{u(t)} = e^{-t\partial_{xxx}} \overline{u_0}$.

\begin{proposition}[Strichartz estimates]\label{free-strichartz} Let $u(t) = e^{-t\partial_{xxx}} u_0$.  Then
\begin{align*}
\| u \|_{L^4_t L^\infty_x} &\lesssim \| u_0 \|_{\dot H^{-1/4}_x}\\
\| u \|_{L^6_{t,x}} &\lesssim \| u_0 \|_{\dot H^{-1/6}_x}\\
\| u \|_{L^\infty_t L^2_x \cap L^8_{t,x} 
 \cap L^6_t L^\infty_x} 
&\lesssim \| u_0 \|_{L^2_x}
\end{align*}
whenever the right-hand sides are finite.
\end{proposition}

\begin{proof} See \cite{kpv}.
\end{proof}

\begin{remark} The free Strichartz estimates (and some basic Littlewood-Paley decomposition) already show that if $u_0 \in \dot H^{-1/6}_x$, then $u$ is locally in $L^4_t L^\infty_x$, and in particular $u^4$ is locally integrable.  This will help ensure that there will be no difficulty interpreting \eqref{eq:gkdv} in the sense of distributions for the regularities under consideration.
\end{remark}

The above Strichartz estimates do not allow us to place a free solution in $L^4_{t,x}(\R \times \R)$, and indeed such control is not available for any non-zero $u_0$.  However, if one assumes some additional frequency separation (for instance, by inserting a suitable bilinear Fourier multiplier)
then one can place \emph{products} of free solutions in $L^2_{t,x}$:

\begin{proposition}[Bilinear Strichartz estimate]\label{bilprop}
 Let $u(t) := e^{-t\partial_{xxx}} u_0$ and $v(t) := e^{-t\partial_{xxx}} v_0$.  Let $m: \R \times \R \to \C$ be any function such that 
\begin{equation}\label{mbound}
 |m(\xi_1,\xi_2)| \lesssim |\xi_1+\xi_2|^{1/2} |\xi_1-\xi_2|^{1/2}.
 \end{equation}
Then we have
\begin{equation}\label{bil}
\begin{split}
\biggl\| \int_\R \int_\R m(\xi_1,\xi_2) 
\hat u(t,\xi_1) e^{ix\xi_1} &\hat v(t,\xi_2) e^{ix\xi_2}\ d \xi_1 d\xi_2 \biggr\|_{L^2_{t,x}}\\ 
&\lesssim \|u_0 \|_{L^2_x} \|v_0 \|_{L^2_x}.
\end{split}
\end{equation}
\end{proposition}

\begin{proof} See \cite[Lemma 1]{gruenrock}.  The claim there is established for $m(\xi_1,\xi_2)$ exactly equal to
$|\xi_1+\xi_2|^{1/2} |\xi_1-\xi_2|^{1/2}$, but the case for general $m$ obeying \eqref{mbound} follows from an identical argument (and also can be deduced from the special case after several applications of Plancherel's theorem).
\end{proof}

In principle (ignoring for now the issues of what happens to derivatives), we can now place quartic nonlinearities such as $u^4$ in energy spaces such as $L^1_t L^2_x$ in the presence of some frequency magnitude separation, because we can place a quadratic term $u \cdot u$ in $L^2_t L^2_x$ and the other two terms in $L^4_t L^\infty_x$.  This argument does not deal directly with the self-interaction case, when all terms in $u^4$ have the same frequency magnitude (e.g. of frequency close to $\pm N$), but in such cases the quartic nonlinearity will either be highly nonresonant (with spacetime frequency $(\tau,\xi)$ close to either $\pm 2(N^3, N)$ or $\pm 4(N^3, N)$) or have very low frequency (and thus be damped by the $\partial_x$ factor in the nonlinearity), and in either case we shall be able to proceed.  For the latter case we shall use the following estimate.

\begin{proposition}[Quartilinear Strichartz estimate]  For $j=1,2,3,4$, let $u_j(t) := e^{-t\partial_{xxx}} u_{j,0}$.  Then
$$ \| (P_+ u_1) (P_+ u_2) (P_- u_3) (P_- u_4) \|_{L^1_t \dot H^{1/2}_x}
\lesssim \prod_{j=1}^4 \| u_{j,0} \|_{\dot H^{-1/4}_x}.
$$
\end{proposition}

\begin{proof} As is well known, the space $\dot H^{1/2}_x \cap L^\infty_x$ is an algebra.  Thus it will suffice to
prove the estimate
$$ \| (P_+ u_1) (P_- u_3) \|_{L^2_t \dot H^{1/2}_x
\cap L^2_t L^\infty_x}
\lesssim \| u_{1,0} \|_{\dot H^{-1/4}_x} \| u_{3,0} \|_{\dot H^{-1/4}_x}$$
and similarly for $u_2,u_4$.  The $L^2_t L^\infty_x$ estimate follows from the $L^4_t L^\infty_x$ Strichartz estimate\footnote{Note that while the Riesz projections $P_\pm$ do not preserve $L^\infty_x$, they will preserve $\dot H^{-1/4}_x$.  Since $P_\pm$ commutes with $e^{t \partial_{xxx}}$, there is thus no problem here.}, so
we are reduced to showing the $L^2_t \dot H^{1/2}_x$ estimate.  Writing $u := |\nabla|^{-1/4} u_{1}$ and
$v := |\nabla|^{-1/4} u_3$, this becomes
\begin{align*}
\| &\int_{-\infty}^0 \int_0^\infty |\xi_1+\xi_2|^{1/2} |\xi_1|^{1/4} |\xi_2|^{1/4} 
\hat u(t,\xi_1) \hat v(t, \xi_2)\ d\xi_1 d\xi_2 \|_{L^2_{t,x}}\\
&\lesssim \| u(0) \|_{L^2_x} \| v(0) \|_{L^2_x},
\end{align*}
but this follows from Proposition \ref{bilprop}.
\end{proof}

\section{The global iteration space}\label{glob-sec}

We define the dyadic nonlinearity space $\N_0 = \N_0(\R \times \R \to \C)$ as the space
$$ \N_0 := L^1_t L^2_x + \dot X^{0,-1/2,1}.$$
We then define the full nonlinearity space $\dot \N^{-1/6} = \dot \N^{-1/6}(\R \times \R \to \C)$ as
$$ \|F\|_{\dot \N^{-1/6}} := \left(\sum_N (N^{-1/6} \| P_N F \|_{\N_0})^2\right)^{1/2},$$
and similarly define $\N^s = \N^s(\R \times \R \to \C)$ for any $s \in \R$ by
$$ \|F\|_{\N^s} := \left(\sum_N (\langle N \rangle^s \| P_N F \|_{\N_0})^2\right)^{1/2}.$$
We also define the dyadic solution space $\S_0 = \S_0(\R \times \R \to \C)$ by the norm
\begin{equation}\label{uso}
\|u\|_{\S_0} := \| u(0) \|_{L^2_x} + \| u_t + u_{xxx} \|_{\N_0}
\end{equation}
and the full solution space $\dot \S^{-1/6} = \dot \S^{-1/6}(\R \times \R \to \C)$ by the norm
$$ \|u\|_{\dot \S^{-1/6}} := \left(\sum_N (N^{-1/6} \| P_N F \|_{\S_0})^2\right)^{1/2}$$
and similarly define $\S^s = \S^s(\R \times \R \to \C)$ for any $s \in \R$ by
$$ \|u\|_{\S^s} := \left(\sum_N (\langle N \rangle^s \| P_N u \|_{\S_0})^2\right)^{1/2}.$$

The advantage of using the spaces $\N_0$ and $\S_0$ is that a large class of estimates for free solutions automatically
extend to $\S_0$ functions. Indeed, we have the following abstract (and well-known) lemma:

\begin{lemma}[Extension lemma]\label{ext}  Let $Y$ be any spacetime Banach space which obeys the time modulation estimate
\begin{equation}\label{gf}
 \| g(t) F(t,x) \|_Y \leq \|g\|_{L^\infty_t(\R)} \|F(t,x)\|_Y
 \end{equation}
for any $F \in Y$ and $g \in L^\infty_t(\R)$.  Let $T: (f_1,\ldots,f_k) \to T(f_1,\ldots,f_k)$ be a spatial multilinear operator for which one has the estimate
$$ \| T( e^{-t\partial_{xxx}} u_{1,0},\ldots, e^{-t\partial_{xxx}} u_{k,0} ) \|_Y \lesssim \prod_{j=1}^k \| u_{j,0} \|_{L^2_x}$$
for all $u_{1,0},\ldots,u_{k,0} \in L^2_x(\R)$.  Then one also has the estimate
$$ \| T( u_1,\ldots, u_k ) \|_Y \lesssim_k \prod_{j=1}^k \| u_j \|_{\S_0}$$
for all $u_1,\ldots,u_k \in \S_0$.
\end{lemma}

\begin{proof} It suffices to prove the claim when $k=1$, as the claim for general $k$ then follows from induction, freezing one of the functions $u_k$ to view $T$ as a $k-1$-multilinear operator, extending the estimate so that $u_1,\ldots,u_{k-1}$ lie in $\S_0$, and then applying the $k=1$ case to extend $u_k$ to $\S_0$ also.  Thus we have
\begin{equation}\label{tuo}
\| T( e^{-\partial_{xxx}} u_0 ) \|_Y \lesssim \|u_0 \|_{L^2_x(\R)}
\end{equation}
for all $u_0 \in L^2_x$, and it will suffice to show that
$$ \| Tu \|_Y \lesssim 1$$
whenever $\|u\|_{\S_0} \lesssim 1$.

By the Duhamel formula and \eqref{uso} we can write
\begin{align*}
u(t) &= e^{-t\partial_{xxx}} u(0) + \int_\R (1_{t > t'} - 1_{0 > t'}) e^{-(t-t')\partial_{xxx}} F_1(t')\ dt' \\
&\quad + \int_0^t e^{-(t-t')\partial_{xxx}} F_2(t')\ dt' 
\end{align*}
where $\|u_0\|_{L^2_x}$, $\|F_1\|_{L^1_t L^2_x}$, and $\|F_2\|_{\dot X^{0,-1/2,1}}$ are all $O(1)$.   The contribution of the first two terms are acceptable from \eqref{tuo} and Minkowski's inequality, using \eqref{gf} to discard the time cutoffs $1_{t > t'} - 1_{0 > t'}$.  For the final term, we split it further as
\begin{align*}
\int_0^t e^{-(t-t')\partial_{xxx}} F_2(t')\ dt'
&= \frac{1}{2} e^{-t\partial_{xxx}} \int_\R \sgn(t') e^{t'\partial_{xxx}} F_2(t')\ dt'\\
&\quad + \frac{1}{2} \int_\R \sgn(t-t') e^{(t-t') \partial_{xxx}} F_2(t')\ dt'.
\end{align*}
For the first term we can use \eqref{uso}, combined with the standard energy estimate
\begin{equation}\label{energy-x}
\| \int_\R \sgn(t') e^{t'\partial_{xxx}} F_2(t')\ dt'\|_{L^2_x} \lesssim \|F_2\|_{\dot X^{0,-1/2,1}}
\end{equation}
which is easily verified from Plancherel's theorem.  For the second term, we use the standard energy estimate
$$ \| \int_\R \sgn(t-t') e^{(t-t') \partial_{xxx}} F_2(t')\ dt' \|_{\dot X^{0,1/2,1}}
\lesssim \|F_2\|_{\dot X^{0,-1/2,1}} $$
which again follows from Plancherel's theorem, and we are reduced to establishing that
$$ \| T v \|_Y \lesssim \| v \|_{\dot X^{0,1/2,1}}$$
for all $v$.  From the definition of the $\dot X^{0,1/2,1}$ norm and the triangle inequality it thus suffices to show that
$$ \| T v \|_Y \lesssim M^{1/2} \| v \|_{L^2_{t,x}}$$
whenever $M > 0$ and the spacetime Fourier transform $\tilde v$ of $v$ is supported in the region $\{ (\tau,\xi):
|\tau - \xi^3| \sim M \}$.  But we can then expand
$$ v = \int_{\lambda \sim M} e^{i\lambda t} e^{-t\partial_{xxx}} v_\lambda\ d\lambda$$
where $v_\lambda(x)$ is defined via the Fourier transform as
$$ \hat v_\lambda(\xi) := \tilde v( \xi^3 + \lambda, \xi ).$$
Applying Minkowski's inequality, \eqref{tuo}, \eqref{gf} we thus have
$$ \| Tv \|_Y \lesssim \int_{\lambda \sim M} \| v_\lambda \|_{L^2_x(\R)}\ d\lambda
\lesssim M^{1/2} (\int_{\lambda \sim M} \| v_\lambda \|_{L^2_x(\R)}^2\ d\lambda)^{1/2}$$
and the claim then follows from Plancherel's theorem.
\end{proof}

As a consequence of this lemma we can extend all the estimates of the previous section to $\S_0$:

\begin{corollary}[$\S_0$ Strichartz estimates]  Let $m: \R \times \R \to \C$ obey \eqref{mbound} and $N > 0$.  Then
\begin{align}
\| u \|_{L^4_t L^\infty_x} &\lesssim \| |\nabla|^{-1/4} u \|_{\S_0}\label{4s}\\
\| u \|_{L^6_{t,x}} &\lesssim \| u \|_{\dot \S^{-1/6}}\label{6s}\\
\| u \|_{L^\infty_t L^2_x \cap L^8_{t,x} \cap L^6_t L^\infty_x} 
&\lesssim \| u \|_{\S_0} \label{lows}\\
\biggl\| \int_\R \int_\R m(\xi_1,\xi_2) 
\hat u(t,\xi_1) e^{ix\xi_1} \hat v(t,\xi_2) e^{ix\xi_2}\ d \xi_1 d\xi_2 & \biggr \|_{L^2_{t,x}} \nonumber\\
&\lesssim \|u_0 \|_{\S_0} \|v_0 \|_{\S_0}\label{bils}\\
\| (P_+ u_1) (P_+ u_2) (P_- u_3) (P_- u_4) \|_{L^1_t \dot H^{1/2}_x(\R \times \R \to \C)}
&\lesssim \prod_{j=1}^4 \| |\nabla|^{-1/4} u_j \|_{\S_0}.\label{quarts}
\end{align}
\end{corollary}

From \eqref{lows} and Bernstein's inequality we obtain the additional estimate
\begin{equation}\label{bern}
\| P_N u \|_{L^\infty_t L^\infty_x} \lesssim N^{1/2} \| u \|_{\S_0}.
\end{equation}

We also need one further $\S_0$ estimate which is easy and standard, though not quite within the purview of 
Lemma \ref{ext} (due to the failure of \eqref{gf} in this setting):

\begin{proposition}[$X^{s,b}$ estimate]\label{xsb}  We have
$$ \| u \|_{\dot X^{0,1/2,\infty}} \lesssim  \| u \|_{\S_0}.$$
\end{proposition}

\begin{proof}  If $u$ is a free solution, $u = e^{-t\partial_{xxx}} u(0)$, then the claim is clear since $\|u\|_{\S_0}$ controls
$\|u(0)\|_{L^2_x}$, and $\dot X^{0,1/2,\infty}$ contains $L^2_x$ free solutions.  By Duhamel's formula, it then suffices to show that
$$ \left\| \int_0^t e^{-(t-t')\partial_{xxx}} F(t')\ dt' \right\|_{\dot X^{0,1/2,\infty}}
\lesssim 1$$
whenever $F$ has unit norm in either $L^1_t L^2_x$ or $\dot X^{0,-1/2,1}$.
We first observe that
$$ \left\| \int_{-\infty}^0 e^{-(t-t')\partial_{xxx}} F(t')\ dt' \right\|_{\dot X^{0,1/2,\infty}}
\lesssim 1.$$
This is because the expression in the left-hand side is an $L^2_x$ free solution, thanks to \eqref{lows}.  Thus it suffices to show that
$$ \left\| \int_{-\infty}^t e^{-(t-t')\partial_{xxx}} F(t')\ dt' \right\|_{\dot X^{0,1/2,\infty}}
\lesssim 1.$$
When $F$ has unit norm in $\dot X^{0,-1/2,1}$, the claim follows by a direct computation involving the spacetime Fourier transform (cf. \eqref{energy-x}).  If instead $F$ has unit norm in $L^1_t L^2_x$, we see from Minkowski's inequality that it suffices to show that
$$ \left\| 1_{t > t'} e^{-(t-t')\partial_{xxx}} F(t') \right\|_{\dot X^{0,1/2,\infty}}
\lesssim \|F(t')\|_{L^2_x}$$
for each $t' \in \R$.  But this can again be seen by a direct computation involving the spacetime Fourier transform\footnote{To be fully rigorous, one must first smoothly truncate the $t$ variable to a compact set, and then take limits, but we omit the details.}.
\end{proof}

As a consequence we can now establish the main nonlinear estimate.

\begin{proposition}[Nonlinear estimate]\label{nonlinear}  For any $u_1,u_2,u_3,u_4 \in \dot \S^{-1/6}$, we have
$$ \| ( u_1 u_2 u_3 u_4 )_x \|_{\dot \N^{-1/6}} \lesssim \prod_{j=1}^4 \| u_j \|_{\dot \S^{-1/6}}.$$
Similarly, for any $s > -1/6$ we have
$$ \| ( u_1 u_2 u_3 u_4 )_x \|_{\N^s} \lesssim_s \sum_{i=1}^4 \| u_i \|_{\S^s} \prod_{j \neq i} 
\| u_j \|_{\dot \S^{-1/6}}.$$
\end{proposition}

\begin{proof} We normalise $\|u_j\|_{\dot \S^{-1/6}} = 1$ for $j=1,2,3,4$ and reduce to showing that
$$ \sum_N [ N^{-1/6} \| P_N (u^4)_x \|_{\N_0}]^2 \lesssim 1$$
and
$$ \sum_N [ \langle N \rangle^s \| P_N (u^4)_x \|_{\N_0}]^2 \lesssim_s 
\sum_{j=1}^4 \sum_N [\langle N \rangle^s\| P_N u_j \|_{\S_0}^2.$$

We write $c_{N,j} := N^{-1/6} \|P_N u_j\|_{\S_0}$, thus
\begin{equation}\label{cnorm}
\sum_N c_{N,j}^2 = 1.
\end{equation}
We perform the Littlewood-Paley decomposition $u_j = \sum_N u_{N,j}$ where $u_{N,j} := P_N u_j$.  After the triangle inequality and permutation invariance, we reduce to showing that
\begin{equation}\label{nnn}
 \sum_N [ \sum_{N_1 \geq N_2 \geq N_3 \geq N_4} 
N^{-1/6} \| P_N ( \prod_{j=1}^4 u_{N_j,j} )_x \|_{\N_0} ]^2 \lesssim 1
\end{equation}
and
\begin{equation}\label{nnn-2}
 \sum_N [ \sum_{N_1 \geq N_2 \geq N_3 \geq N_4} 
\langle N^s \rangle \| P_N ( \prod_{j=1}^4 u_{N_j,j} )_x \|_{\N_0} ]^2 \lesssim 
\sum_N [\langle N^s \rangle N^{1/6} c_{N,1}]^2.
\end{equation}

We can prove \eqref{nnn} and \eqref{nnn-2} simultaneously by establishing the estimate
\begin{equation}\label{nnn-unif}
N^{-1/6} \| P_N ( \prod_{j=1}^4 u_{N_j,j} )_x \|_{\N_0}
\lesssim \min\left( \left(\frac{N}{N_1}\right)^{-\sigma}, \left(\frac{N}{N_1}\right)^{-\sigma}\right) \left(\frac{N_2}{N_4}\right)^{-\sigma} \prod_{j=1}^4 c_{N_j,j}
\end{equation}
whenever $N_1 \geq N_2 \geq N_3 \geq N_4$,
for some explicit constant $\sigma > 0$ (in fact we will have $\sigma = 1/12$).  Indeed, to show \eqref{nnn} it now suffices
to show that
$$ \sum_N \left[ \sum_{N_1 \geq N_2 \geq N_3 \geq N_4}
\min\left( \left(\frac{N}{N_1}\right)^{-\sigma}, \left(\frac{N_1}{N}\right)^{-\sigma}\right) \left(\frac{N_2}{N_4}\right)^{-\sigma} \prod_{j=1}^4 c_{N_j,j} \right]^2 \lesssim 1.$$
But by bounding $c_{N_3,3}$ crudely by $1$ and using Young's inequality and \eqref{cnorm} we see that
$$ \sum_{N_2 \geq N_3 \geq N_4} (N_2/N_4)^{-\sigma} c_{N_2,2} c_{N_3,3} c_{N_4,4} \lesssim 1$$
and then the claim then follows from another application of Young's inequality and \eqref{cnorm}.  A similar argument gives
\eqref{nnn-2}.

It remains to prove \eqref{nnn-unif}.
Note that we may assume that $N_1 \gtrsim N$ since the expression in the norm vanishes otherwise.
Let us first consider the ``non-self-interaction'' case when $N_1/N_4 > 1.001$ are not adjacent.  Here we estimate the $\N_0$ norm by the $L^1_t L^2_x$ norm, and use the Strichartz estimates \eqref{4s}, \eqref{bils} to obtain
\begin{align*}
\| P_N ( \prod_{j=1}^4 u_{N_j,j}  )_x \|_{\N_0} 
&\leq \| P_N ( \prod_{j=1}^4 u_{N_j,j} )_x \|_{L^1_t L^2_x} \\
&\lesssim N \| \prod_{j=1}^4 u_{N_j,j} \|_{L^1_t L^2_x} \\
&\lesssim N \| u_{N_1,1} u_{N_4,4} \|_{L^2_t L^2_x} \| u_{N_2,2} \|_{L^4_t L^\infty_x} \|u_{N_3,3} \|_{L^4_t L^\infty_x} \\
&\lesssim N N_1^{-1} \| u_{N_1,1} \|_{\S_0} \| u_{N_4,4} \|_{\S_0} N_2^{-1/4} \| u_{N_2,2} \|_{\S_0} N_3^{-1/4} \| u_{N_3,3} \|_{\S_0} \\
&= N N_1^{-5/6} N_2^{-1/12} N_3^{-1/12} N_4^{1/6} c_{N_1,1} c_{N_2,2} c_{N_3,3} c_{N_4,4}
\end{align*}
which implies \eqref{nnn-unif} with $\sigma = 1/12$ (bounding $N_3^{-1/12}$ by $N_4^{-1/12}$).

Now we consider the case when $N_1,N_4$ are adjacent, thus 
$$ N_4 \leq N_3 \leq N_2 \leq N_1 \leq 1.001 N_4.$$
Suppose first that one of the $u_{N_j}$ has spacetime Fourier transform vanishing on the set
$$\{ (\tau,\xi): |\tau - N_j^2 \sgn(\xi)| \leq 0.01 N_j^2 \}.$$
From Proposition \ref{xsb} we then have
$$ \|u_{N_j}\|_{L^2_t L^2_x} \lesssim N_j^{1/6} N_j^{-3/2} c_{N_j}$$
while from \eqref{lows} we have
$$ \| u_{N_k} \|_{L^6_t L^\infty_x} \lesssim N_k^{1/6} c_{N_k}$$
for the other three values of $k$.  By H\"older, Bernstein, and the comparability of the $N_j$ we then have
$$ \| P_N ( u_{N_1,1} u_{N_2,2} u_{N_3,3} u_{N_4,4} )_x \|_{L^1_t L^2_x}
\lesssim N N_1^{-5/6} c_{N_1,1} c_{N_2,2} c_{N_3,3} c_{N_4,4}$$
and \eqref{nnn-unif} follows in this case.  Thus by smooth Fourier decomposition, we may assume that
$u_{N_j}$ has Fourier transform supported on the set
$$\{ (\tau,\xi): |\tau - N_j^2 \sgn(\xi)| \leq 0.02 N_j^2 \}$$
for each $j=1,2,3,4$.

Next, we use Riesz transforms in space to split $u_{N_j,j} = P_+ u_{N_j,j} + P_- u_{N_j,j}$ for $j=1,\ldots,4$, thus giving sixteen terms in the product $u_{N_1,1} \ldots u_{N_4,4}$.  Consider any term in which there are more $P_+$ than $P_-$ or vice versa.
Then some elementary algebra (and the fact that $N_1, \ldots, N_4$ differ by at most a factor of $1.001$) 
shows that this component of $u_{N_1,1} \ldots u_{N_4,4}$ has spacetime Fourier transform supported at a distance at least $\gtrsim N_1^3$ from the cubic $\tau = \xi^3$.  Thus we may bound this contribution to \eqref{nnn-unif} by
$$
N^{-1/6} N N_1^{-3/2} \| P_{\pm} u_{N_1,1} P_{\pm} u_{N_2,2} P_{\pm} u_{N_3,3} P_{\pm} u_{N_4,4} \|_{L^2_t L^2_x}$$
where the four signs $\pm$ need not be equal.
From \eqref{lows} (and the boundedness of the Riesz transforms) we have
$$ \| P_\pm u_{N_j,j} \|_{L^8_t L^8_x} \lesssim N_j^{1/6} c_{N_j,j}$$
for $j=1,\ldots,4$, so we can bound the above expression by
$$ N^{5/6} N_1^{-5/6} c_{N_1,1} c_{N_2,2} c_{N_3,3} c_{N_4,4}$$
which is acceptable.  Thus we only need to consider the case when there are two $P_+$'s and two $P_-$'s.  After some relabeling, it suffices to show that
$$
N^{5/6} \| P_N ( (P_+ u_1) (P_+ u_2) (P_- u_3) (P_- u_4) ) \|_{L^1_t L^2_x} 
\lesssim \left(\frac{N_1}{N}\right)^{-\sigma} 
\prod_{j=1}^4 N_j^{1/6} \| u_j \|_{\S_0} 
$$
whenever $u_j$ have Fourier support in the region $|\xi| \sim N_1$.  
Applying \eqref{quarts} we see that the left-hand side is
$$ N^{1/3} N_1^{-1/4} \|u_1 \|_{\S_0} N_1^{-1/4} \| u_2 \|_{\S_0}
N_1^{-1/4} \| u_3 \|_{\S_0} N_1^{-1/4} \| u_4 \|_{\S_0}]^2$$
and the claim follows.
\end{proof}

Using this estimate and a standard iteration argument, we can thus conclude that if $\| u_0 \|_{\dot H^{-1/6}_x}$ is sufficiently small, then there exists a unique solution $u \in \S_0$ with small norm to \eqref{eq:gkdv} (interpreted of course in the Duhamel sense).  Further standard techniques\footnote{See for instance \cite{kpv} for very analogous arguments in $L^2_x$ for the quintic gKdV equation.} then show that $u$ varies analytically (as a function from $\dot H^{-1/6}_x$ to $\S_0$) with respect to the data $u_0$, again using the small norm assumption.  Also one can use standard continuity arguments to strengthen the uniqueness claim slightly\footnote{It may be possible to strengthen the uniqueness claim further.  For instance, for the KdV equation, weak solutions in the $L^\infty_t L^2_x$ class are shown to be unique in \cite{zhou}.  We will not pursue this issue here.}, so that the solution $u$ is unique among all solutions in $\S_0$ (not necessarily of small norm); we omit the details.  From the second estimate in Proposition \ref{nonlinear} and standard arguments we can establish persistence of regularity, or more precisely if $u_0$ has small $\dot H^{-1/6}_x$ norm and is also in $H^s_x$ for some $s > -1/6$, then $u(t)$ will lie in $H^s_x$ for all time, and in fact we have the estimate
$$ \| u \|_{L^\infty_t H^s_x} \lesssim_s \| u_0 \|_{H^s_x}.$$
We omit the details as they are very standard (see e.g. \cite{kpv} for very similar arguments).
These results in particular imply that the solutions constructed here are the unique strong $\dot S^{-1/6}$ limit of classical (smooth) solutions, and are thus compatible with the solutions constructed previously at higher regularities in \cite{kpv} or \cite{gruenrock}.

One can also show scattering for these solutions in $\dot H^{-1/6}_x$ norm, thus there exists $u_+ \in \dot H^{-1/6}_x$ such that
$u(t) - e^{-t\partial_{xxx}} u_+$ converges in $\dot H^{-1/6}_x$ as $t \to +\infty$.  Indeed, from the Duhamel formula
$$ u(t) = e^{-t\partial_{xxx}} u_0 - e^{-t\partial_{xxx}} \int_0^t e^{t' \partial_{xxx}} F(t')\ dt'$$
where $F := (u^4)_x$, it suffices to show that the integral $\int_0^{+\infty} e^{t' \partial_{xxx}} F(t')\ dt'$ is
conditionally convergent in $\dot H^{-1/6}_x$.  But from Proposition \ref{nonlinear} we know
that $F \in \dot \N^{-1/6}$, and so the claim would follow from

\begin{lemma}\label{scatterlemma} Let $F \in \dot \N^{-1/6}$.  Then $\int_0^{+\infty} e^{t' \partial_{xxx}} F(t')\ dt'$ is conditionally convergent in $\dot H^{-1/6}_x$.
\end{lemma}

\begin{proof} The claim is clear when $F$ is a test function.  Such functions are easily verified to be dense in $\dot \N^{-1/6}$.  From the energy estimate (in \eqref{lows}, say) we know that the $L^\infty_t \dot \N^{-1/6}$ norm of
$\int_0^t e^{t' \partial_{xxx}} F(t')\ dt'$ is controlled by the $\N^{-1/6}$ norm of $F$. The claim then follows from standard limiting arguments.
\end{proof}

If one furthermore has $u_0 \in H^s_x$ in addition to being small in $\dot H^{-1/6}_x$, then an adaptation of the persistence of regularity argument also shows that $u(t) - e^{-t\partial_{xxx}} u_+$ converges in $H^s_x$ norm; again, we omit the details.

\section{Large data local well-posedness}\label{largedata}

In the preceding section we showed that for any initial data $u_0$ with sufficiently small $\dot H^{-1/6}_x$ norm, that there was a unique global solution $u$ to the Cauchy problem \eqref{eq:gkdv} with small $\dot \S^{-1/6}$ norm, with the solution depending continuously on the data.  To complete the proof of Theorem \ref{main} we need to also establish local wellposedness for large $\dot H^{-1/6}_x$.  The precise claim is as follows.

\begin{theorem}[Critical local existence for large data]\label{cle}  For any $R > 0$ there exists a constant $\eps > 0$ such that the following claim is true: whenever $u_0 \in \dot H^{-1/6}_x$ and $N > 0$ is such that $\|u_0\|_{\dot H^{-1/6}_x} \leq R$ and $\| P_{>N} u_0 \|_{\dot H^{-1/6}_x} \leq \eps$, then there exists a unique solution $u \in \dot \S^{-1/6}( [- \eps N^{-3}, \eps N^{-3}] \times \R )$ to \eqref{eq:gkdv} on the time interval $[- \eps N^{-3}, \eps N^{-3}]$ with the bounds
$$\| u \|_{\dot \S^{-1/6}( [- \eps N^{-3}, \eps N^{-3}] \times \R )} \lesssim R$$
and
$$ \| P_{>N} u \|_{\dot \S^{-1/6}( [- \eps N^{-3}, \eps N^{-3}] \times \R )} \lesssim \eps.$$
Furthermore, for data $u_0$ restricted to the above class, the map $u_0 \mapsto u$ from data to solution is Lipschitz continuous from $\dot H^{-1/6}_x$ to $\dot \S^{-1/6}( [- \eps N^{-3}, \eps N^{-3}] \times \R )$.
\end{theorem}

We will only sketch the proof here, as it is a fairly standard modification of the global existence arguments.  First we may 
use a scaling argument to normalise $N=1$.  The main task is then to show that the nonlinear map $u \mapsto \Phi(u)$ defined
as
$$ \Phi(u)(t) := e^{-t\partial_{xxx}} u_0 + \int_0^t e^{-(t-t')\partial_{xxx}} \partial_x( u^4(t') )\ dt'$$
is a contraction on the space
$$ \{ u: \| u \|_{\dot \S^{-1/6}( [- \eps, \eps] \times \R)} \lesssim R; \quad
\| P_{>1} u \|_{\dot \S^{-1/6}( [- \eps, \eps] \times \R)} \lesssim \eps \}$$
with suitable choices of implied constants, and with a metric given by the norm
\begin{equation}\label{ir}
 \frac{1}{R} \| u \|_{\dot \S^{-1/6}( [- \eps, \eps] \times \R)} + \frac{1}{\eps} \| P_{>1} u \|_{\dot \S^{-1/6}( [- \eps, \eps] \times \R)}.
\end{equation}
To prove this contraction property, one needs a variant of Proposition \ref{nonlinear}, but with the global
norm $\dot \S^{-1/6}(\R \times \R)$ replaced by the variant \eqref{ir}, and similarly for
$\dot \N^{-1/6}(\R \times \R)$.  The only issue is to ensure that the final bound gains at least a fractional power of $\eps$ to counteract any factors of $R$ which appear.  If at least two frequencies in the quartic nonlinearity are greater than $1$, then this gain of $\eps$ is automatic from the arguments used to prove Proposition \ref{nonlinear}.  The only new feature arises when at most one of the factors has frequency greater than one.  However, in this case we can take advantage of the localisation of the time interval to $[-\eps,\eps]$.  For instance, from \eqref{bern} and a H\"older in time we have
$$ \| P_N u \|_{L^4_t L^\infty_x([-\eps,\eps] \times \R)} \lesssim \eps^{1/4} N^{1/2} \| u \|_{\S_0([-\eps,\eps] \times \R )}.$$
This estimate is superior to \eqref{4s} when $N \leq 1$, as it gains a fractional power of $\eps$ (and also gains some powers of $N$).  If one uses this estimate as a substitute for \eqref{4s} in Proposition \ref{nonlinear}, we can already gain the desired power of $\eps$ in almost all cases, except for the case in which all four frequencies are equal (or adjacent) 
and less than $1$.  But in this case one can argue by many means, for instance one can use $L^\infty_t L^\infty_x$ and
$L^\infty_t L^2_x$ estimates to place the nonlinearity in (say) $L^\infty_t L^{3/2}_x$, which after a H\"older in time
places one in $L^1_t L^{3/2}_x$ (gaining the desired power of $\eps$), and then Sobolev embedding places one in $L^1_t \dot H^{-1/6}_x$, which will suffice.  We omit the details.  This concludes the proof of Theorem \ref{main}.

\begin{remark} The local existence theorem also gives a blowup criterion for large $\dot H^{-1/6}_x(\R \to \C)$ solutions to \eqref{eq:gkdv}; if the solution cannot be continued past some time $T_*$, this means that either
$\limsup_{t \to T_*^-} \| u(t) \|_{\dot H^{-1/6}_x} = \infty$, or that
$\liminf_{N \to \infty} \liminf_{t \to T_*^-} \| P_{>N} u(t) \|_{\dot H^{-1/6}_x} > 0$.  In other words, either we have norm blowup or else the $\dot H^{-1/6}_x$ concentrates in increasingly higher frequencies. It is likely that a refinement of the arguments here would produce a more usable blowup criterion; for instance, a reasonable conjecture would be that blowup (or more generally, failure of scattering) would only occur if the $L^6_{t,x}$ norm of $u$ was infinite.  We will however not pursue this matter here (and in any case, we already have global wellposedness of real-valued solutions of regularity $L^2_x(\R)$ or higher).
\end{remark}

\begin{remark} It is likely that one could combine the methods here (or the simpler methods in \cite{gruenrock}) with the ``$I$-method'' (as used for instance in \cite{Ckstt}) to push the large (real-valued) data global existence results for to 
regularities below $L^2_x(\R)$.  However, it is unlikely that these methods would get arbitrarily close to the scaling regularity $s=-1/6$.  We will not pursue this issue here.
\end{remark}

\section{Proof of Theorem \ref{main2}}\label{main2-sec}

Now we begin the proof of Theorem \ref{main2}.  Let the notation and assumptions be as in Theorem \ref{asymptotic} and
Theorem \ref{main2}.  We write
\begin{equation}\label{rt}
R(t) := \frac{1}{\lambda^{2/3}(t)} Q\left( \frac{x-x(t)}{\lambda(t)} \right)
\end{equation}
thus $u = R + w$.  Thus $w$ solves the forced gKdV equation
\begin{equation}\label{weq}
w_t + w_{xxx} + (w^4)_x = E
\end{equation}
where the error term $E$ is defined by
$$ E := ( R^4 + w^4 - (R+w)^4 )_x - (R_t + R_{xxx} + (R^4)_x ).$$
Now we estimate the error.  In fact it is exponentially localised in $L^2_{t,x}$ to the soliton trajectory
$\{ (t,x(t)): t \in \R \}$.

\begin{lemma}[Error estimate] We have
\begin{equation}\label{eee}
\int_\R \int_\R |E(t,x)|^2 e^{|x-x(t)|}\ dx dt \lesssim \eps^2.
\end{equation}
\end{lemma}

\begin{proof} Let us first control the contribution of the $( R^4 + w^4 - (R+w)^4 )_x$ term.  A simple application of
the product rule gives
$$ | ( R^4 + w^4 - (R+w)^4 )_x | \lesssim |R_x| (|w|^3 + R^2 |w|) + |w_x| (|w|^2 R + R^3).$$
From the definition of $R$, $Q$, and \eqref{lambda} we thus see that (if $\eps$ is sufficiently small)
$$ | ( R^4 + w^4 - (R+w)^4 )_x | \lesssim e^{-0.9 |x-x(t)|} ( |w| + |w_x| ) ( 1 + |w|^2 ).$$
From Sobolev embedding and \eqref{lh1} we have the pointwise estimate
$$ |w(t,x)| \lesssim \|w\|_{L^\infty_t H^1_x(\R \times \R \to \C)} \lesssim \eps$$
and so 
$$ | ( R^4 + w^4 - (R+w)^4 )_x |^2 e^{|x-x(t)|} \lesssim e^{-0.8|x-x(t)|} ( |w|^2 + |w_x|^2 )$$
and the claim follows from \eqref{wox}.

It remains to control the contribution of $R_t + R_{xxx} + (R^4)_x$.  Direct computation
(using \eqref{rt}, \eqref{soliton-eq}) yields
$$ R_t + R_{xxx} + (R^4)_x
= -\frac{2}{3} \frac{\lambda'(t)}{\lambda(t)} R - \frac{\lambda'(t)}{\lambda(t)} (x-x(t)) R_x
- \frac{x'(t) - \lambda(t)^{-2}}{\lambda(t)} R_x $$
and hence by \eqref{lambda}, \eqref{lambda-rate} and the exponential decay of $Q$ (and hence $R$) we have the pointwise estimate
$$ |\partial_t R + \partial_{xxx} R + \partial_x(R^4)|(t,x) \lesssim e^{-0.9|x-x(t)|}
(\int_\R |w(t,y)|^2 e^{-\frac{1}{2} |y-x(t)|}\ dy)^{1/2}.$$
Using \eqref{wox} we thus see that this contribution is also acceptable, and we are done.
\end{proof}

At time zero we have $w(0) = (u(0)-Q) + (Q-R)$, and thus by the hypotheses on $u(0)-Q$ and \eqref{lambda} we easily see that
\begin{equation}\label{woo}
 \|w(0)\|_{H^1_x(\R)} + \|w(0)\|_{\dot H^{-1/6}_x(\R)} \lesssim \eps.
\end{equation}

The strategy is now to modify the proof of Theorem \ref{main} to ``solve'' the forced gKdV equation \eqref{weq}.  Unfortunately, the forcing term $E$ does not quite lie in any of the $\N$ family of spaces, and so $w$ will not lie in the $\S$
family.  Thus we shall have to expand these spaces in somewhat complicated ways, dependent on the function $x(t)$.  
For any $h \in \R$, let $\Gamma_h \subset \R \times \R$ be the set
$$ \Gamma_h := \{ (t,x): |x - x(t) - h| \leq 1 \}.$$
This set thus tracks the trajectory of the soliton, shifted in space by a fixed shift $h$. 

\begin{definition}[Modified nonlinearity space]  We define $\N_*$ to be the Banach space generated by atoms $F$ from one of the following two types:
\begin{itemize}
\item (Standard atoms) We have $\|F\|_{\dot \N^{-1/6}(\R \times \R \to \C) \cap \N^1(\R \times \R \to \C)} \leq 1$.
\item (Exotic atoms) There exists an $h \in \R$ such that $F$ is supported on $\Gamma_h$ such that $\|F\|_{L^2_t L^2_x(\Gamma_h \to \C)} \leq 1$.
\end{itemize}
\end{definition}

By covering spacetime by regions $\Gamma_h$, where $h$ ranges over integers, we see from \eqref{eee}
that
\begin{equation}\label{eman}
\|E\|_{\N_*} \lesssim \eps.
\end{equation}
Thus one can view $\N_*$ as basically the minimal extension to $\N^1 \cap \dot \N^{-1/6}$ that can accomodate forcing terms such as $E$.  The exotic atoms are \emph{almost} standard - for instance, one can place them in $\dot X^{s,-1/2,\infty}$ rather than $\dot X^{s,-1/2,1}$ for $-1/6 \leq s \leq 1$ - but for the critical problem it seems unfortunately necessary to treat these exotic atoms as a separate component of the nonlinearity space.  Note that even though Theorem \ref{main2} is concerned only with real-valued functions, we permit $\N_*$ to contain complex-valued functions; this is a technicality arising from our use of the Fourier transform.

For the space $\S_*$, the heuristic motivation for the space is as follows.  From the Duhamel formula we see that forcing terms $F$ which lie in $L^1_t L^2_x$ spaces lead to solutions $u$ which can be viewed as superpositions of free solutions truncated
to half-spaces such as $\{ (t,x): t \geq t' \}$.  The exotic atoms generate forcing terms which can be viewed as lying in a ``curvilinear'' variant of the $L^1_t L^2_x$, and heuristically one expects those forcing terms to generate solutions which act like (superpositions of) free solutions, truncated to curvilinear half-spaces such as $\{ (t,x): x \leq x(t) + h \}$.  Actually for technical reasons it is more convenient to work with smooth truncations to such spaces, thus introducing an additional error term which is localised to one of the regions $\Gamma_h$.

More precisely, we shall need to fix a smooth cutoff function $\eta: \R \to \R$ with $\eta(x) = 1$ for $x \leq -1$, $\eta(x) = 0$ for $x \geq +1$, and $\eta$ smoothly interpolated in between; it is also convenient to impose the symmetry constraint
\begin{equation}\label{eta-sym}
\eta(x) + \eta(-x) \equiv 1.
\end{equation}
We then define $\S_*$ as follows:

\begin{definition}[Modified solution space]  We define $\S_*$ to be the Banach space generated by atoms $u$ from one of the following three types:
\begin{itemize}
\item (Standard atoms) We have $\| u \|_{\dot \S^{-1/6}(\R \times \R \to \C) \cap \S^1(\R \times \R \to \C)} \leq 1$.
\item (Semi-standard atoms) We have $u(t,x) = \tilde u(t,x) \eta(x - x(t) - h)$ for some $h \in \R$, where $\tilde u$ is a standard atom of $\S_*$.
\item (Exotic atoms) There exists an $h \in \R$ such that $u$ is supported on $\Gamma_h$, and $\|u\|_{L^2_t H^1_x(\Gamma_h \to \C) + L^\infty_t H^1_x(\R \times \R \to \C)} \leq 1$.
\end{itemize}
\end{definition}

\begin{remark} The standard atoms are in fact redundant, being the limit of semi-standard atoms as $h \to +\infty$.  However we retain them for expository purposes to emphasise that $\S_*$ is an extension of $\dot \S^{-1/6} \cap \S^1$.
\end{remark}

In the next section we shall establish the following two basic estimates on $\S_*$ and $\N_*$:

\begin{proposition}[Modified energy estimate]\label{energy-mod} For any $u: \R \times \R \to \C$, we have
$$ \| u \|_{\S_*} \lesssim \| u(0) \|_{\dot H^{-1/6}_x \cap \dot H^1_x}
+ \| u_t + u_{xxx} \|_{\N_*}.$$
\end{proposition}

\begin{proposition}[Modified nonlinear estimate]\label{nonlinear-mod} For any $u_1,u_2,u_3,u_4 \in \S_*$, we have
$$ \| ( u_1 u_2 u_3 u_4 )_x \|_{\N_*} \lesssim \| u_1 \|_{\S_*} \| u_2 \|_{\S_*} \| u_3 \|_{\S_*} \| u_4 \|_{\S_*}.$$
\end{proposition}

From the modified energy estimate and \eqref{weq}, \eqref{woo}, \eqref{eman} we have
$$ \| w \|_{\S_*} \lesssim \eps + \| ( w^4 )_x \|_{\N_*}$$
and thus by the modified nonlinear estimate
$$ \| w \|_{\S_*} \lesssim \eps + \| w \|_{\S_*}^4.$$
A standard continuity argument then gives (for $\eps$ sufficiently small)
$$ \| w \|_{\S_*} \lesssim \eps$$
and then from the modified nonlinear estimate again
$$ \| w_t + w_{xxx} \|_{\N_*} \lesssim \eps.$$
Now, from Proposition \ref{energy-mod} (and the observation that $\S_*$ controls $L^\infty_t \dot H^{-1/6}_x \cap L^\infty_t H^1_x$) we easily modify the proof of Lemma \ref{scatterlemma} to obtain

\begin{lemma}\label{scatterlemma2} Let $F \in \N_*$.  Then $\int_0^\infty e^{t' \partial_{xxx}} F(t')\ dt'$ is conditionally convergent in $\dot H^{-1/6}_x \cap H^1_x$.
\end{lemma}

From this and the Duhamel formula we thus obtain Theorem \ref{main2}.

\begin{remark} One can modify the above arguments to obtain the higher regularity estimates
$$ \| w \|_{L^\infty_t H^s_x(\R \times \R \to \R)} \lesssim_s \|w(0) \|_{H^s_x(\R \to \R)}$$
for all $s > 1$; we omit the details.
\end{remark}

\section{Proof of estimates}\label{estimates-sec}

In this section we allow our functions to be complex-valued unless otherwise stated.

To finish the proof of Theorem \ref{main2} we need to establish Proposition \ref{energy-mod} and Proposition \ref{nonlinear-mod}.  Both estimates will rely on the following global-in-time, local-in-space Kato smoothing effect
(a nonlinear variant of which formed a crucial component of the arguments in \cite{martel2}).

\begin{proposition}[Time-global Kato local smoothing effect]\label{nrg}  Let $h, h' \in \R$.
\begin{itemize}
\item For any $u_0$, we have
\begin{equation}\label{utx}
 \| e^{-t\partial_{xxx}} u_0 \|_{L^2_t H^1_x(\Gamma_h)} \lesssim \|u_0\|_{L^2_x}.
\end{equation}
\item For any $F$ supported on $\Gamma_{h'}$, we have
\begin{equation}\label{ftx}
 \| \int_{-\infty}^t e^{-(t-t') \partial_{xxx}} F(t')\ dt' \|_{L^\infty_t H^1_x(\R \times \R) + L^2_t H^2_x(\Gamma_h)}
\lesssim \|F\|_{L^2_{t,x}(\Gamma_{h'})}.
\end{equation}
\end{itemize}
\end{proposition}

\begin{proof}
We begin with \eqref{utx}.  We may translate $h=0$ and take $u_0$ to be real-valued.  Write $u := e^{-t\partial_{xxx}} u_0$.  Since $u_t + u_{xxx} = 0$, we observe the mass transport identity
$$ (u^2)_t + (u^2)_{xxx} = 3 \partial_x( u_x^2 ).$$
We multiply this by a function $\psi(x - x(t))$, where $\psi: \R \to \R$ is to be determined shortly, and integrate in $x$ to obtain
\begin{align*}
 \partial_t \int_\R \psi(x-x(t)) u(t,x)^2\ dx &=-  \int_\R (\psi' - \psi''')(x-x(t)) u(t,x)^2\ dx\\
&\quad - 3 \int_\R \psi'(x-x(t)) u_x(t,x)^2\ dx.
\end{align*}
If we now choose $\psi(x) := \tanh(x/100)$, we see that $\psi' - \psi'''$ and $\psi'$ are non-negative, and strictly positive for $-1 < x < 1$, thus
$$ \partial_t \int_\R \tanh((x-x(t))/100) u(t,x)^2\ dx \leq - c\int_{|x-x(t)| < 1} u(t,x)^2 + u_x(t,x)^2\ dx$$
for some absolute constant $c > 0$.
Integrating this in $t$ and using mass conservation we obtain \eqref{utx}.

From \eqref{utx} and a duality argument we can easily estimate the $L^\infty_t H^1_x$ term of \eqref{ftx}, so we concentrate on the $L^2_t H^2_x$ term.  By a limiting argument we may take $F$ to be smooth and compactly supported in spacetime; we may also take $F$ to be real-valued and normalise the $L^2_{t,x}$ norm of $F$ to be $1$.  We also normalise $h=0$.  Let
$u(t) := \int_{-\infty}^t e^{-(t-t') \partial_{xxx}} F(t')\ dt'$, thus $u_t+u_{xxx} = F$ and $u(-\infty) = 0$.  From the preceding discussion we already have
$$ \| u \|_{L^\infty_t H^1_x(\R \times \R)} \lesssim 1.$$
Now $u$ obeys the mass transport identity
$$ (u^2)_t + (u^2)_{xxx} = 3 \partial_x( u_x^2 ) + 2uF$$
and hence
\begin{align*}
\partial_t \int_\R \psi(x-x(t)) u(t,x)^2\ dx &=-  \int_\R (\psi' - \psi''')(x-x(t)) u(t,x)^2\ dx\\
&\quad - 3 \int_\R \psi'(x-x(t)) u_x(t,x)^2\ dx\\
&\quad + 2\int_\R \psi(x-x(t)) u(t,x) F(t,x)\ dx.
\end{align*}
Let us take $\psi(x) := 1 + \tanh(x/100)$ if $h' \leq 0$, and $\psi(x) := 1 - \tanh(x/100)$ if $h' > 0$.  Then a simple application of Cauchy-Schwarz, exploiting the support of $F$, gives
$$ 2\int_\R \psi(x-x(t)) u(t,x) F(t,x)\ dx \leq \frac{1}{2} \int_\R (\psi' - \psi''')(x-x(t)) u(t,x)^2\ dx
+ O( \int_\R F(t,x)^2\ dx )$$
and hence
$$ \partial_t \int_\R \psi(x-x(t)) u(t,x)^2\ dx \leq - c\int_{|x-x(t)| < 1} u(t,x)^2 + u_x(t,x)^2\ dx
+ O( \int_\R F(t,x)^2\ dx )$$
for some $c > 0$.  Integrating this in time we obtain the $L^2_t H^1_x$ bounds.  To get the $L^2_t H^2_x$ bound, we apply the variant mass transport identity
$$ (u_x^2)_t + (u_x^2)_{xxx} = 3 \partial_x( u_{xx}^2 ) + 2
u_x F_x$$
and argue similarly to before, integrating by parts to move the derivative off of $F$; we omit the details.
\end{proof}

\subsection{Proof of Proposition \ref{energy-mod}.}

If $u$ is a free solution (so $u_t + u_{xxx}=0$), then the claim follows 
from \eqref{uso}, as $u$ is just a constant multiple of a standard atom.  So we may assume that
$u(0) = 0$, and reduce to showing that
$$ \| \int_0^t e^{-(t-t')\partial_{xxx}} F(t')\ dt' \|_{\S_*}
\lesssim \|F\|_{\N_*}.$$
We may assume that $F$ is one of the atoms of $\N_*$.  If it is a standard atom then the claim again follows from
\eqref{uso}, so we may assume $F$ is an exotic atom.  After a spatial translation, it suffices to show that
$$ \| \int_0^t e^{-(t-t')\partial_{xxx}} F(t')\ dt' \|_{\S_*}
\lesssim \|F\|_{L^2_{t,x}(\Gamma_0)}$$
whenever $F$ is supported on $\Gamma_0$.  

The intuition, based on the heuristic that the Airy fundamental solution propagates to the left, is that
the forward propagator $\int_\R 1_{t' > t} e^{-(t-t')\partial_{xxx}} F(t')\ dt'$ should be mostly concentrated on the
region $\{ (t,x): x \leq x(t) \}$, while the backward propagator $\int_\R 1_{t' < t} e^{-(t-t')\partial_{xxx}} F(t')\ dt'$
should be concentrated in the region $\{ (t,x): x \geq x(t) \}$. Thus the retardation cutoff $1_{t' > t}$ is morally 
interchangeable with the (time-dependent) spatial cutoff $1_{x \leq x(t)}$, the key point being that the latter does not depend on $t'$.

We now make these heuristics precise. Observe that
\begin{align*}
 \| \int_{-\infty}^0 e^{-(t-t')\partial_{xxx}} F(t')\ dt' \|_{\S_*}
&\lesssim
 \| \int_{-\infty}^0 e^{-(t-t')\partial_{xxx}} F(t')\ dt' \|_{\dot \S^{-1/6} \cap \S^1} \\
&\lesssim \| \int_{-\infty}^0 e^{t' \partial_{xxx'}} F(t')\ dt' \|_{\dot H^{-1/6}_x \cap H^1_x} \\
&\lesssim \|F\|_{L^2_{t,x}(\Gamma_0)}
\end{align*}
thanks to \eqref{uso} and Proposition \ref{nrg}.  Thus it suffices to show that
$$ \| \int_{-\infty}^\infty 1_{t' < t} e^{-(t-t')\partial_{xxx}} F(t')\ dt' \|_{\S_*}
\lesssim \|F\|_{L^2_{t,x}(\Gamma_0)}.$$
A similar argument shows that
$$ \| \int_{-\infty}^{+\infty} e^{-(t-t')\partial_{xxx}} F(t')\ dt' \|_{\dot \S^{-1/6} \cap \S^1}
\lesssim \|F\|_{L^2_{t,x}(\Gamma_0)}$$
and hence
$$ \| \int_{-\infty}^{+\infty} e^{-(t-t')\partial_{xxx}} F(t')\ dt' \eta( x - x(t)) \|_{\S_*}
\lesssim \|F\|_{L^2_{t,x}(\Gamma_0)}.$$
Thus it suffices to show that
$$ \| \int_{-\infty}^\infty [1_{t' < t} - \eta(x - x(t)) ] e^{-(t-t')\partial_{xxx}} F(t')\ dt' \|_{\S_*}
\lesssim \|F\|_{L^2_{t,x}(\Gamma_0)}.$$
We smoothly decompose $\R \times \R$ into the regions $\Gamma_h$ for $h \in \Z$ and rely entirely on exotic atoms to
estimate the left-hand side somewhat crudely as
$$ \sum_{j=0}^1 \sum_{h \in \Z} \| \partial_x^j \int_{-\infty}^\infty [1_{t' < t} - \eta(x - x(t)) ] e^{-(t-t')\partial_{xxx}} F(t')\ dt' \|_{L^2_{t,x}( \Gamma_h )}.$$
For the terms where $-10 \leq h \leq 10$, we may use Proposition \ref{nrg} to obtain satisfactory bounds, treating the $1_{t' < t}$ and $\eta(x-x(t))$ term separately.  Now let us consider terms in which $h > 10$.  Here the $\eta(x-x(t))$ term disappears.  We expand the propagator $e^{-(t-t')\partial_{xxx}}$ as
$$ \partial_x^j [e^{-(t-t')\partial_{xxx}} F(t')](x) = \int_\R [\int_\R e^{i(t-t')\xi^3} e^{i(x-x')\xi} \xi^j\ d\xi] F(t',x')\ dx'.$$
where the inner integral is interpreted in a principal value sense.  Now as $F$ is supported in $\Gamma_0$, we have
$|x' - x(t')| \leq 1$.  On $\Gamma_h$, we have $|x-x(t)-h| \leq 1$, while from \eqref{lambda} we have $x(t) - x(t') \geq \frac{1}{2} (t-t')$, and thus
$$x-x' \geq h  + \frac{1}{2} (t-t') + 2.$$
Standard stationary phase integration by parts arguments (or standard
bounds on the Airy function) then establish the bound
$$ | \int_\R e^{i(t-t')\xi^3} e^{i(x-x')\xi} \xi^j\ d\xi | \lesssim h^{-100} \langle t - t' \rangle^{-100}$$
(say) for $j=0,1$, $(t',x') \in \Gamma_0$, $(t,x) \in \Gamma_h$, and $t > t'$; from this it is easy to see that the contribution of the $h > 10$ terms will be acceptable.

Now we consider the terms when $h < -10$.   Here we use the symmetry \eqref{eta-sym} to rewrite
$$ 1_{t' < t} - \eta(x-x(t)) = - [ 1_{-t' < -t} - \eta( (-x) - (-x(t)) ) ].$$
If we then use the time reversal symmetry $(t,x) \mapsto (-t,-x)$, which maps $x(t)$ to $-x(-t)$, we see that this case now follows from the preceding case.  This concludes the proof of Proposition \ref{energy-mod}.

\subsection{Proof of Proposition \ref{nonlinear-mod}.}

It suffices to show that
$$ \| ( u_1 u_2 u_3 u_4 )_x \|_{\N_*} \lesssim 1$$
whenever $u_1,u_2,u_3,u_4$ are atoms of $\S_*$, possibly of different types.

Let us first see what happens when at least one of the atoms, say $u_1$, is an exotic atom, say on $\Gamma_h$.  Then
$( u_1 u_2 u_3 u_4 )_x$ lies in $\Gamma_h$ and so we may use the exotic $\N_*$ atoms to estimate
$$ \| ( u_1 u_2 u_3 u_4 )_x \|_{\N_*} \lesssim \|( u_1 u_2 u_3 u_4 )_x \|_{L^2_t L^2_x(\Gamma_h)}.$$
Since $H^1_x$ is closed under multiplication, we may estimate the right-hand side by
$$ \lesssim \| u_1 \|_{L^2_t H^1_x(\Gamma_h)} \prod_{j=2}^4 \| u_j \|_{L^\infty_t H^1_x(\R \times \R)}$$
which is easily seen to be bounded by $O(1)$ as desired.  

Thus we may now reduce to the case when all the $u_j$ are either standard or semi-standard.  In fact, as the standard atoms are the limit of semi-standard atoms, we may take each of the $u_j$ to be semi-standard, thus for each $1 \leq j \leq 4$ we have $u_j(t,x) = v_j(t,x) \eta(x - x(t) - h_j)$ for some $v_j$ with 
\begin{equation}\label{jbound}
\| v_j \|_{\dot \S^{-1/6} + \S^1} \leq 1
\end{equation}
and $h_j \in \R$.  Our task is thus to show that
$$ \| ( v_1 v_2 v_3 v_4 \phi )_x \|_{\N_*} \lesssim 1$$
where $\phi := \prod_{j=1}^4 \eta(x-x(t)-h_j)$.
We split the left-hand side into the paraproducts
\begin{equation}\label{para-1}
\| \sum_N P_N ( v_1 v_2 v_3 v_4 P_{< N/100} \phi )_x \|_{\N_*}
\end{equation}
and
\begin{equation}\label{para-2}
 \| \sum_N P_N ( v_1 v_2 v_3 v_4 P_{\geq N/100} \phi )_x \|_{\N_*}.
\end{equation}
Consider the first paraproduct.  We estimate the $\N_*$ norm using standard atoms and reduce to showing that
$$ ( \sum_N [(N^{-1/6} + \langle N \rangle) 
\| P_N ( v_1 v_2 v_3 v_4 P_{< N/100} \phi )_x \|_{\N_0}]^2)^{1/2} \lesssim 1.$$
This can be achieved by direct modification of the arguments in Proposition \ref{nonlinear}; the point is that
the new factor $P_{< N/100} \phi$ is bounded, and sufficiently localised in spatial frequency, that it will not disrupt any
of the estimates involved in the proof of Proposition \ref{nonlinear}.  We omit the details, and turn instead to the second paraproduct.  Here we use the triangle inequality and reduce to showing that
$$
\sum_N \| P_N ( v_1 v_2 v_3 v_4 P_{\geq N/100} \phi )_x \|_{\N_*} \lesssim 1.$$
The operator $P_N \partial_x$ is a spatial convolution operator with a kernel of $L^1_x$ norm $O(N)$, so we reduce to
showing
$$
\sum_N N \| v_1 v_2 v_3 v_4 P_{\geq N/100} \phi \|_{\N_*} \lesssim 1.$$
We use exotic atoms and cover $\R \times \R$ into regions $\Gamma_h$ for $h \in \Z$ and reduce to showing
$$
\sum_{h \in \Z} \sum_N N \| v_1 v_2 v_3 v_4 P_{\geq N/100} \phi \|_{L^2_{t,x}(\Gamma_h)} \lesssim 1.$$
From Proposition \ref{nrg} and Lemma \ref{ext} we have
$$ \| v_j \|_{L^\infty_t H^1_x(\R \times \R)} + \| v_j \|_{L^2_t H^1_x(\Gamma_h)} \lesssim 1$$
for $j=1,2,3,4$.
Since $H^1_x$ is closed under multiplication, we conclude in particular that
$$ \|v_1 v_2 v_3 v_4 \|_{L^2_t L^2_x(\Gamma_h)} \lesssim 1$$
and so by H\"older's inequality it suffices to show that
$$
\sum_{h \in \Z} \sum_N N \| P_{\geq N/100} \phi \|_{L^\infty_{t,x}(\Gamma_h)} \lesssim 1.$$
However, direct computation (and integration by parts) shows the pointwise estimate
\begin{align*}
|P_{\geq N/100} \phi(t,x)| &\lesssim \sum_{j=1}^4 \langle N \rangle^{-100} \langle N (x - x(t)-h_j) \rangle^{-100} \\
&\lesssim \langle N \rangle^{-100} \sum_{j=1}^4 \langle N (h-h_j) \rangle^{-100}
\end{align*}
and the claim follows.  This concludes the proof of Proposition \ref{nonlinear-mod}.


\begin{thebibliography}{10}

\bibitem{bej}
I. Bejenaru, \emph{On Schr\"odinger Maps}, preprint.

\bibitem{cct}
M. Christ, J. Colliander, T. Tao, \emph{Asymptotics, frequency modulation, and low regularity ill-posedness for canonical defocussing equations}, Amer. J. Math. \textbf{125} (2003), 1235--1293.

\bibitem{ck}
M. Christ, A. Kiselev, \emph{Maximal operators associated to filtrations}, J. Funct. Anal. \textbf{179} (2001).

\bibitem{Ckstt}
J. Colliander, M. Keel, G. Staffilani, H. Takaoka, T. Tao, \emph{Sharp global well-posedness results for periodic and non-periodic KdV and modified KdV on $R$ and $T$}, J. Amer. Math. Soc. \textbf{16} (2003), 705--749.

\bibitem{cote}
R. C\^ote, \emph{Large data wave operator for the generalized Korteweg-de Vries equations}, Differential and Integral Equations, \textbf{19} (2006), no. 2, 163--188.

\bibitem{cote2}
R. C\^ote, \emph{Construction of solutions to the subcritical gKdV equations with a given asymptotic behaviour}, preprint.

\bibitem{cote3}
R. C\^ote, \emph{Construction of solutions to the $L^2$-critical KdV equation with a given asymptotic behaviour}, preprint.

\bibitem{gruenrock}
A. Gr\"unrock, \emph{A bilinear Airy-estimate with application to gKdV-3}, Differential Integral Equations \textbf{18} (2005), no. 12, 1333--1339.

\bibitem{ik}
A. Ionescu, C. Kenig, \emph{Low-regularity Schr\"{o}dinger maps, II: global well-posedness in dimensions $d\geq 3$}, preprint.

\bibitem{kpv}
C. Kenig, G. Ponce, L. Vega, \emph{Wellposedness and scattering results for the generalized Korteweg-de Vries equation via the contraction principle}, Comm. Pure Appl. Math. \textbf{46} (1993), 527--560.

\bibitem{martel}
Y. Martel, F. Merle, \emph{Asymptotic stability of solitons for subcritical generalized KdV equations}, Arch. Ration. Mech. Anal. \textbf{157} (2001), no. 3, 219--254. 

\bibitem{martel2}
Y. Martel, F. Merle, \emph{Asymptotic stability of solitons of the subcritical gKdV equations revisited}, Nonlinearity \textbf{18} (2005), no. 1, 55--80.

\bibitem{hideo}
H. Takaoka, Y. Tsutsumi, \emph{Well-posedness of the Cauchy problem for the modified KdV equation with periodic boundary condition}. Int. Math. Res. Not. \textbf{56} (2004), 3009--3040.

\bibitem{tao:wavemap2}
T. Tao, \emph{Global regularity of wave maps II.  Small energy in two dimensions}, submitted, Comm. Math. Phys.

\bibitem{tataru:wave2}
D. Tataru, \emph{On global existence and scattering for the wave maps equation},
Amer. J. Math. \textbf{123} (2001), no. 1, 37--77.

\bibitem{zhou}
Y. Zhou, \emph{Uniqueness of weak solution of the KdV equation}, Int. Math. Res. Not. \textbf{6} (1997), 271--283.

\end{thebibliography}
\end{document}